\documentclass[12pt]{article}

\usepackage{graphicx}

\def\ra{\rightarrow}

 \def\HollowBox #1#2{{\dimen0=#1 \advance\dimen0 by -#2       
       \dimen1=#1 \advance\dimen1 by #2                       
        \vrule height #1 depth #2 width #2                    
        \vrule height 0pt depth #2 width #1                   
        \llap{\vrule height #1 depth -\dimen0 width \dimen1}%
       \hskip -#2                                             
       \vrule height #1 depth #2 width #2}}                   
 \def\BoxOpTwo{\mathord{\HollowBox{6pt}{.4pt}}\;}             

\def\endpf{\hfill $\BoxOpTwo$}

\font\teneufm=eufm10
\font\seveneufm=eufm7
\font\fiveeufm=eufm5
\newfam\eufmfam
\textfont\eufmfam=\teneufm
\scriptfont\eufmfam=\seveneufm
\scriptscriptfont\eufmfam=\fiveeufm

\newfam\msbfam
\font\tenmsb=msbm10  scaled \magstep1 \textfont\msbfam=\tenmsb
\font\sevenmsb=msbm7 scaled \magstep1 \scriptfont\msbfam=\sevenmsb
\font\fivemsb=msbm5  scaled \magstep1 \scriptscriptfont\msbfam=\fivemsb
\def\Bbb{\fam\msbfam \tenmsb}

\def\RR{{\Bbb R}}

\newfam\msbbfam
\font\tenmsbb=msbm10  scaled \magstep3 \textfont\msbbfam=\tenmsbb
\font\sevenmsbb=msbm7 scaled \magstep3 \scriptfont\msbbfam=\sevenmsbb
\font\fivemsbb=msbm5  scaled \magstep3 \scriptscriptfont\msbbfam=\fivemsbb
\def\Bbbb{\fam\msbbfam \tenmsbb}

\def\RRR{{\Bbbb R}}

\newtheorem{theorem}{Theorem}
\newtheorem{corollary}[theorem]{Corollary}
\newtheorem{proposition}[theorem]{Proposition}
\newtheorem{lemma}[theorem]{Lemma}

\makeindex

\begin{document}

\begin{center}
\huge \bf Higher Dimensional Conundra
\end{center}

\begin{center}
\large Steven G. Krantz\footnote{We are happy to thank
the American Institute of Mathematics for its hospitality
and support during this work.}
\end{center}
\vspace*{.15in}

\begin{quote}
{\bf Abstract:}  In recent years, especially in
the subject of harmonic analysis, there has been interest
in geometric phenomena of $\RR^N$ as $N \rightarrow +\infty$. 
In the present paper we examine several specific geometric
phenomena in Euclidean space and calculate the asymptotics
as the dimension gets large.
\end{quote}

\setcounter{section}{-1}

\section{Introduction}

Typically when we do geometry we concentrate on a specific venue in 
a particular space.  Often the context is {\it Euclidean space}, and
often the work is done in $\RR^2$ or $\RR^3$.  But in modern
work there are many aspects of analysis that are linked to concrete
aspects of geometry.  And there is often interest in rendering the
ideas in Hilbert space or some other infinite dimensional setting.
Thus one wants to see how the finite-dimensional result in $\RR^N$ changes
as $N \ra +\infty$.

In the present paper we study some particular aspects of the geometry of
$\RR^N$ and their asymptotic behavior as $N \ra \infty$.  We choose these
particular examples because the results are surprising or especially
interesting.   One may hope that they will lead to further studies.

\section{Volume in \boldmath $\RRR^N$}

Let us begin by calculating the volume of the unit ball in $\RR^N$ and
the surface area of its bounding unit sphere.   We let $\Omega_N$ denote
the former and $\omega_{N-1}$ denote the latter.  In addition, we let
$\Gamma(x)$ be the celebrated Gamma function of L. Euler.   It is 
a helpful intuition (which is literally true when $x$ is an integer) that
$\Gamma(x) \approx (x - 1)!$.  We shall also use Stirling's formula which
says that
$$
k! \approx k^k \cdot e^{-k} \cdot \sqrt{2\pi k} 
$$
or, more generally,
$$
\Gamma(x) \approx (x-1)^{x-1} e^{-(x-1)} \sqrt{2\pi(x-1)} 
$$
for $x \in \RR$, $x > 0$.

\begin{lemma}   \sl 
We have that
$$ 
\int_{\RR^N} e^{- \pi\|{\bf x}\|^2} \, d{\bf x} = 1 . 
$$
\end{lemma}
{\bf Proof:}  The case $N = 1$ is familiar from calculus.  We write
$$
S = \int_\RR e^{-\pi t^2} \, dt
$$
hence
\begin{eqnarray*}
S^2 & = & \int_\RR e^{-\pi x^2} \, dx \int_\RR e^{-\pi y^2} \, dy \\ 
    & = & \mathop{\int\!\!\!\int}_{\RR^2} e^{-\pi \|{\bf x}\|^2} d{\bf x}  \\
    & \stackrel{\rm (polar\ coordinates)}{=} & \int_0^\infty \int_{\|{\bf x}\| = 1} e^{-\pi r^2} r \, ds({\bf x}) dr  \\
    & = & \omega_1 \frac{1}{2\pi} e^{-\pi r^2} \biggr |_0^\infty  \\
    & = & \frac{\omega_1}{2\pi} \\\
    & = & 1 
\end{eqnarray*}
hence $S = 1$.

For the $N-$dimensional case, write
$$ 
\int_{\RR^N} e^{-\pi |x|^2} dx = \int_{\RR} e^{- \pi x_1^2} dx_1 \cdots \int_{\RR} e^{- \pi x_N^2} dx_N 
$$
and apply the one-dimensional result. 
\endpf 
\medskip \\

Let $\sigma$ be the unique rotationally invariant area measure on $S_{N-1} = \partial B_N$.

\begin{lemma}    \sl
We have
$$  
\omega_{N-1} = \frac{2\pi^{N/2}}{\Gamma(N/2)} \, . 
$$
\end{lemma}
{\bf Proof:} \ \  Introducing polar coordinates we have
$$ 
1 = \int_{\RR^N} e^{- \pi|x|^2} dx = 
                                    \int_{S^{N-1}} d\sigma \int_0^\infty e^{-\pi r^2} r^{N-1} dr 
$$
or 
$$ 
\frac{1}{\omega_{N-1}} = \int_0^\infty 
                              e^{-\pi r^2} r^{N} \frac{dr}{r} . 
$$
Letting $s = r^2$ in this last integral and doing some obvious manipulations
yields the result.  
\endpf 
\medskip \\

\begin{corollary} \sl
The volume of the unit ball in $\RR^N$ is
$$
\Omega_N = \frac{\omega_{N-1}}{N} = \frac{2\pi^{N/2}}{\Gamma(N/2) \cdot N} \, .
$$ 
\end{corollary}
{\bf Proof:}  We calculate that
$$
\Omega_N  = \mathop{\int\!\!\!\int}_B 1 \, dV(x) \stackrel{\rm (polar\ coordinates)}{=} 
\int_0^1 \int_{\|{\bf x}\| = 1} 1 \cdot r^{N-1} \, d\sigma({\bf x}) dr = \omega_{N-1} \cdot \frac{r^N}{N} \biggr |_0^1 =
  \frac{\omega_{N-1}}{N} \, .
$$
That completes the proof.
\endpf
\medskip \\

Now the first nontrivial fact that we wish to observe about the volume of the Euclidean
unit ball in $N$-space is that that volume tends to 0 at $N \ra \infty$.  More formally,

\begin{proposition} \sl
We have the limit
$$
\lim_{N \ra +\infty} \Omega(N) = 0 \, .
$$
\end{proposition}
{\bf Proof:}   We calculate that
\begin{eqnarray*}
\hbox{\rm (Volume\ of\ Unit\ Ball)} & = & \frac{2 \pi^{N/2}}{\Gamma(N/2) \cdot N} \\
     & \approx & \frac{2 \pi^{N/2}}{((N-2)/2)^{(N-2)/2} e^{-(N-2)/2}\sqrt{2\pi [(N-2)/2]} \cdot N}  \\
     & \approx & \frac{(2\pi e)^{N/2} \cdot 2}{N^{(N-1)/2} \cdot \sqrt{\pi} \cdot N} \\
     & \approx & \frac{(2\pi e)^{N/2} \cdot 2}{N^{(N+1)/2} \cdot \sqrt{\pi}} \\
     & \approx & \left ( \frac{2\pi e}{N} \right )^{N/2} \cdot \frac{1}{\sqrt{N}} \cdot \frac{2}{\sqrt{\pi}} \, .
\end{eqnarray*}
This expression clearly tends to 0 as $N \ra +\infty$.
\endpf 
\medskip \\

In fact we can actually say something about the {\it rate} at which the volume of the
ball tends to zero.  We have

\begin{proposition} \sl
We have the estimate
$$
0 \leq \Omega_N \leq 2 \cdot \frac{20^{N/2}}{N^{(N+1)/2}} \, .
$$
\end{proposition}
{\bf Proof:}  Follows by inspection of the last line of the proof of Proposition 4.
\endpf 
\medskip \\

In fact something more is true about the volumes of balls in high-dimensional
Euclidean space.

\begin{proposition} \sl
Let $R > 0$ be fixed.  Then
$$
\lim_{N \ra +\infty} \hbox{Vol}(B(0,R)) = 0 \, .
$$
In other words, the volume of the ball of radius $R$ tends to 0.
\end{proposition}
{\bf Proof:}  From the formula for the volume of the unit ball we have that
$$
\lim_{N \ra +\infty} \hbox{Vol}(B(0,R)) =
  \lim_{N \ra +\infty} \left ( \frac{2\pi e R^2}{N} \right )^{N/2} \cdot \frac{1}{\sqrt{N}} \cdot \frac{2}{\sqrt{\pi}} \, .
$$
This expression clearly tends to 0 as $N \ra +\infty$.
\endpf
\medskip \\

We leave the proof of the next result as an exercise for the reader; simply examine
the formula for $\omega_{N-1}$:

\begin{proposition} \sl
Let $R > 0$.  Then the surface area of the sphere of radius $R$ in $\RR^N$ tends to
0 as $N \ra +\infty$.
\end{proposition}

The following very simple but remarkable fact comes up in considerations
of spherical summation of Fourier series.

\begin{proposition} \sl
As $N \ra +\infty$, the volume of the unit ball in $\RR^N$ is concentrated more
and more out near the boundary sphere.  More precisely, let $\delta > 0$.  Then
$$
\lim_{N \ra +\infty} \frac{\hbox{\rm volume}(B(0,1) \setminus B(0,1 - \delta))}{\hbox{\rm volume}(B(0,1)} = 1 \, .
$$
\end{proposition}
{\bf Proof:}  We have
\begin{eqnarray*}
\lim_{N \ra +\infty} \frac{\hbox{\rm volume}(B(0,1) \setminus B(0,1 - \delta))}{\hbox{\rm volume}(B(0,1)} & = &
     \lim_{N \ra + \infty} \frac{[1 - (1 - \delta)^N] \cdot [2\pi^{N/2}]/[\Gamma(N/2) \cdot N]}{[2\pi^{N/2}]/[\Gamma(N/2) \cdot N]} \\
      & = & \lim_{N \ra +\infty} 1 - (1 - \delta)^N \\
      & = & 1 \, .
\end{eqnarray*}
That is the desired conclusion.
\endpf
\medskip \\

\section{A Case of Leakage}

The title of this section gives away the punchline of the example.  Or so it may seem to some.

Consider at first a square box of side two with sides parallel to the coordinate axes in the Euclidean plane.
We may inscribe in this box four discs of diameter 1, as shown in Figure 1.  These discs will be
called {\it primary discs}.  Once those four discs
are inscribed, we may inscribe a small, shaded disc in the middle as shown in Figure 2.
We set
$$
{\cal R}_2 = \frac{\hbox{area\ of\ shaded\ disc}}{\hbox{area\ of\ large\ box}} \, .
$$

    \begin{figure}
    \centering
      \includegraphics[height=2.65in, width=2.75in]{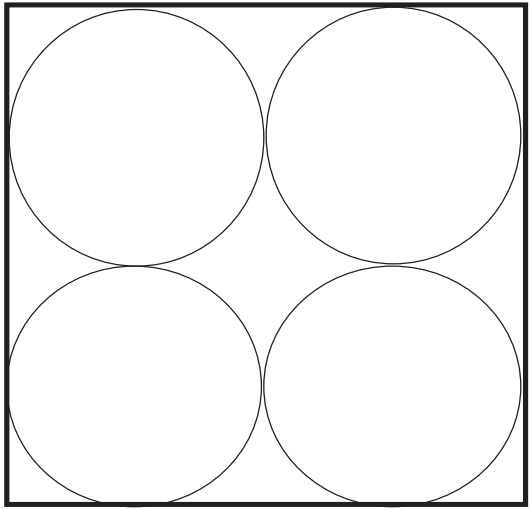}
      \caption{The configuration in dimension 2.}
    \end{figure}

    \begin{figure}
    \centering
      \includegraphics[height=2.65in, width=2.75in]{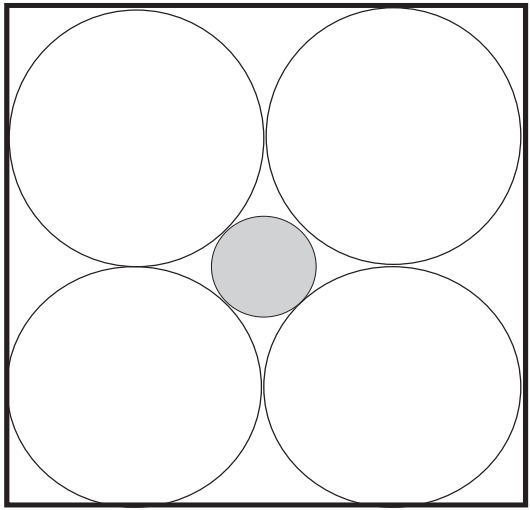}
      \caption{The shaded disc in dimension 2.}
    \end{figure}

The same construction may be performed in Euclidean dimension 3.  Examine Figure 3.  It suggests a rectangular
parallelepiped with all sides equal to 2, and 8 unit balls inscribed inside in a canonical fashion.  These
eight primary balls determine a unique inscribed shaded ball in the center.  We set
$$
{\cal R}_3 = \frac{\hbox{volume\ of\ shaded\ ball}}{\hbox{volume\ of\ large\ box}} \, .
$$

    \begin{figure}
    \centering
      \includegraphics[height=2.65in, width=3.2in]{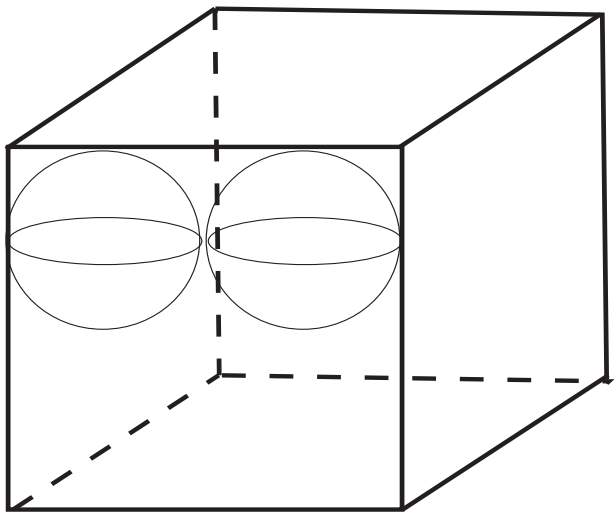}
      \caption{The configuration in dimension 3.}
    \end{figure}

A similar construction may be performed in any dimension $N \geq 2$, with $2^N$ balls inscribed
in a rectangular box of side 2.  The ratio ${\cal R}_N$ is then calculated in just the same
way.  The question is then
\begin{quote}
What is the limit \ \ $\lim_{N \ra\infty} {\cal R}_N$ \ as \ $N \ra +\infty$?
\end{quote}

It is natural to suppose, and most people do suppose, and that this limit (assuming it exists) is
between 0 and 1.  All other things being equal, it is likely equal to either 0 or 1.  Thus
it comes as something of a surprise that this limit is in fact equal to $+\infty$.  Let us now
enunciate this result and prove it.

\begin{proposition} \sl
The limit
$$
\lim_{N \ra +\infty} {\cal R}_N = +\infty \, .
$$
\end{proposition}

Of course this result is counter-intuitive, because we all instinctively believe that the shaded ball, in
any dimension, is contained inside the big box.  Such is not the case.  We are being fooled by
the 2-dimensional situation depicted in Figure 1.  In that special situation, any of the two adjacent
primary discs actually touch in such a way as to trap the shaded disc in a particular convex subregion of the big box
(see Figure 4).  So certainly it must be that ${\cal R}_2 < 1$.  But such is not the case in higher
dimensions.  There is actually a gap on each side of the box through which the shaded ball can leak.
And indeed it does.

    \begin{figure}
    \centering
      \includegraphics[height=2.65in, width=2.75in]{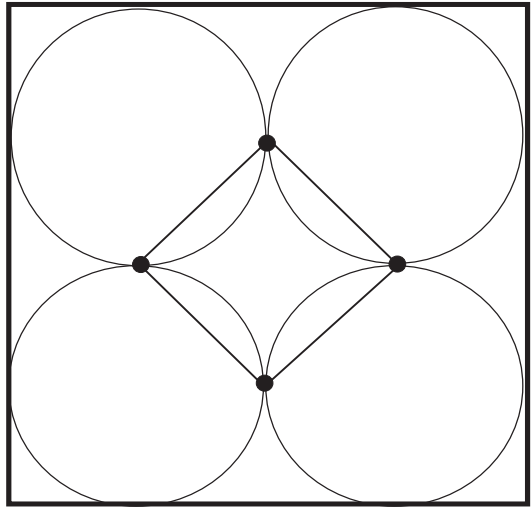}
      \caption{The disc trapped in dimension 2.}
    \end{figure}

This is what we shall now show.  First we shall perform the calculation of ${\cal R}_N$ for each
$N$ and confirm that the expression tends to $+\infty$ as $N \ra +\infty$.  Then we shall
calculate the first dimension in which the shaded ball actually leaks out of the box.  
\medskip \\

\noindent {\bf Proof of the Proposition:} 
Notice that the center of one of the primary balls is at the point $(1,1, \dots, 1)$. 
It is a simple matter to calculate that a boundary point of this ball that is
nearest to the center of the box is located at $P^* \equiv (1 - 1/\sqrt{N}, 1 - 1/\sqrt{N}, \dots, 1 - 1/\sqrt{N}$.
Since the shaded ball will osculate the primary ball at that point, we see that
the shaded ball has center the origin and radius equal to
$$
\hbox{dist}(0, P^*) = \sqrt{ (1 - 1/\sqrt{N})^2,  1 - 1/\sqrt{N})^2, \dots, 1 - 1/\sqrt{N})^2 } = \sqrt{N + 1 - 2\sqrt{N}} \, .
$$

Thus we see that the volume of the shaded ball is
$$
[N + 1 - 2\sqrt{N}]^{N/2} \cdot \Omega_N \, .
$$
The ratio ${\cal R}_N$ is then
$$
{\cal R}_N = \frac{[N + 1 - 2\sqrt{N}]^{N/2} \cdot \Omega_N}{2^N} \, .
$$

Now we may simplify this last expression to
$$
{\cal R}_N = \frac{2 \cdot \pi^{N/2}}{\Gamma(N/2) \cdot N} \cdot \frac{[N + 1 - 2\sqrt{N}]^{N/2}}{2^N} \, .
$$
After some simplification we find that
$$
{\cal R}_N = \frac{2 (\pi/4)^{N/2} \cdot [N + 1 - 2\sqrt{N}]^{N/2}}{\Gamma(N/2) \cdot N} \, .
$$
By Stirling's formula, this last expression is approximately equal to
\begin{eqnarray*}
\lefteqn{\frac{2 \cdot (\pi/4)^{N/2} (N + 1 - 2\sqrt{N})^{N/2}}{N} \cdot \left ( \frac{N-2}{2} \right )^{(2-N)/2} \cdot e^{(N-2)/2} \cdot
  \frac{1}{\sqrt{\pi(N-2)}}} \\
& = &
  \frac{2}{N} \left [ \frac{\pi}{4} (N + 1 - 2\sqrt{N}) \cdot \frac{2}{N-2} \cdot e \right ]^{N/2} \cdot
    \frac{N-2}{2} \cdot \frac{1}{\sqrt{\pi(N-2)}} \cdot \frac{1}{e} \, .
\end{eqnarray*}
After some manipulation, we finally find that
\begin{eqnarray*}
\lim_{N \ra \infty} {\cal R}_N & = & \lim_{N \ra \infty} \frac{2}{N} \left ( \frac{\pi e}{2} \right )^{N/2} \cdot
  \left ( \frac{N+1}{N-2} \right )^{N/2} \cdot \frac{N-2}{2} \cdot \frac{1}{\sqrt{\pi(N-2)}} \cdot \frac{1}{e} \\
  & = & \lim_{N \ra \infty} \frac{2}{N} \left ( \frac{\pi e}{2} \right )^{N/2} \left ( 1 + \frac{3}{N-2} \right )^{N/2} \cdot
    \frac{N-2}{2} \cdot \frac{1}{\sqrt{\pi(N-2)}} \cdot \frac{1}{e} \, .
\end{eqnarray*}
Now, in the limit, we may replace expressions like $N-2$ by $N$.  And we may reparametrize $N$
as $3N$.  The result is
\begin{eqnarray*}
\lefteqn{\lim_{N \ra \infty} \frac{2}{N} \left ( \frac{\pi e}{2} \right )^{N/2} \cdot \left (
    1 + \frac{3}{N} \right )^{N/2} \cdot \frac{N-2}{2} \cdot \frac{1}{\sqrt{\pi(N-2)}} \cdot \frac{1}{e}}  \\
     & = &
    \lim_{N \ra \infty} \frac{2}{N} \left ( \frac{\pi e}{2} \right )^{N/2} \cdot
      \left ( 1 + \frac{1}{N} \right )^{3N/2} \cdot \frac{N-2}{2} \cdot \frac{1}{\sqrt{\pi(N-2)}} \cdot \frac{1}{e} \\
        & = &
   \lim_{N \ra \infty} \frac{2}{N} \cdot \left ( \frac{\pi e}{2} \right )^{N/2} \left [ \left ( 1 + \frac{1}{N} \right )^N  \right ]^{3/2} \cdot
        \frac{N-2}{2} \cdot \frac{1}{\sqrt{\pi(N-2)}} \cdot \frac{1}{e} \, .
\end{eqnarray*}

What we see now is that this last equals
$$
\lim_{N \ra \infty} \left ( \frac{\pi e}{2} \right )^{N/2} \cdot \frac{1}{\sqrt{N}} \cdot \frac{e^{1/2}}{\sqrt{\pi}} \, .
$$
Plainly, because $\pi e/2 > 4$, this limit is $+\infty$.  That proves the result.
\endpf
\medskip \\

And now we turn to the question of when the shaded ball starts to leak out of the big box.  This
is in fact easy to analyze.  We need only determine when the radius of the shaded ball exceeds 1.
First notice that the radius of the shaded ball is monotone increasing in $N$.  Now we need
to solve
$$
\sqrt{N + 1 - 2\sqrt{N}} > 1 \, .
$$
This is a simple algebra problem, and the solution is $N > 4$.  Thus, beginning in dimension 5,
the shaded ball will ``leak out of'' the large box.

It may be noted that Richard W. Cottle has made a study of mathematical phenomena that change
(in the manner of a catastrophe---see [ZEE]) between dimensions 4 and dimensions 5.  The results
 may be found in [COT].

\section{Centroids}

This final section of the paper will be more like an invitation to further
exploration.  We cannot include all the details of the calculations, as they
are too recondite and complex.  Yet the topic is very much in the spirit of
the theme of this paper, and we cannot resist including a few pointers
to this new and interesting work (for which see [KRA1] and [KRMP]).

The inspiration for this work is the following somewhat surprising observation.
Let $T$ be a triangle in the plane (see Figure 5).  There are three ways
to calculate the centroid of this figure:  {\bf (i)}  average the vertices,
{\bf (ii)} average the edges, or {\bf (iii)} average the 2-dimensional solid
figure.  And the question is:  are these three versions of the centroid the same?
The answer is that {\bf (i)} and {\bf (iii)} are {\it always the same}.   Generically
{\bf (ii)} is different.  In fact the three versions of the centroid coincide if
and only if the triangle is equilateral [KRMP].

    \begin{figure}
    \centering
      \includegraphics[height=2.65in, width=2in]{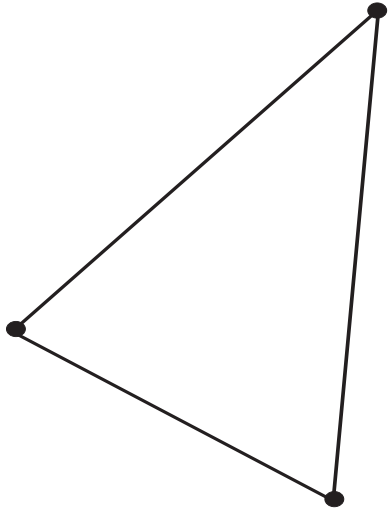}
      \caption{Centroids for a triangle.}
    \end{figure}

We used this fact as a springing-off point to investigate analogous questions in higher
dimensions.  Consider the simplex {\bf S} in $\RR^N$ that is the convex hull of the points
$0 = (0,0, \dots, 0)$, $(1, 0, \dots, 0)$, $(0, 1, 0, \dots, 0)$, \dots, $(0,0, \dots,0, 1)$.
Refer to Figure 6.
Such an $N$-dimensional geometric figure comes equipped with $(N+1)$ notions of centroid:
one can average the vertices (or 1-dimensional skeleton) ${\cal S}_0$, or one can average over the 1-dimensional skeleton
${\cal S}_1$, or
one can average over the two-dimensional skeleton ${\cal S}_2$A, or \dots one can average over
the $(N-1)$-dimensional skeleton ${\cal S}_{N-1}$, or one can average over the $N$-dimensional solid ${\cal S}_N$.
There results the centroids ${\cal C}_{0,N}$, ${\cal C}_{1, N}$, \dots, ${\cal C}_{N,N}$.
And the question is:  Are these different notions of centroid all the same?
And here is the somewhat surprising answer:
\begin{quote}
In dimensions 2 through 12 (for the ambient space), the skeletons ${\cal S}_0$ and ${\cal S}_N$ 
for the simplex {\bf S} have the same centroid.
In those same dimensions, the skeletons ${\cal S}_1$, ${\cal S}_2$, \dots, ${\cal S}_{N-1}$ 
all have different
centroids, and the centroids all differ from the common centroid for ${\cal S}_0$ and ${\cal S}_N$.
But in dimension 13 things are different.  In fact in that dimension the skeletons ${\cal S}_3$ and
${\cal S}_8$ have the same centroid.
\end{quote}

    \begin{figure}
    \centering
      \includegraphics[height=2.65in, width=2in]{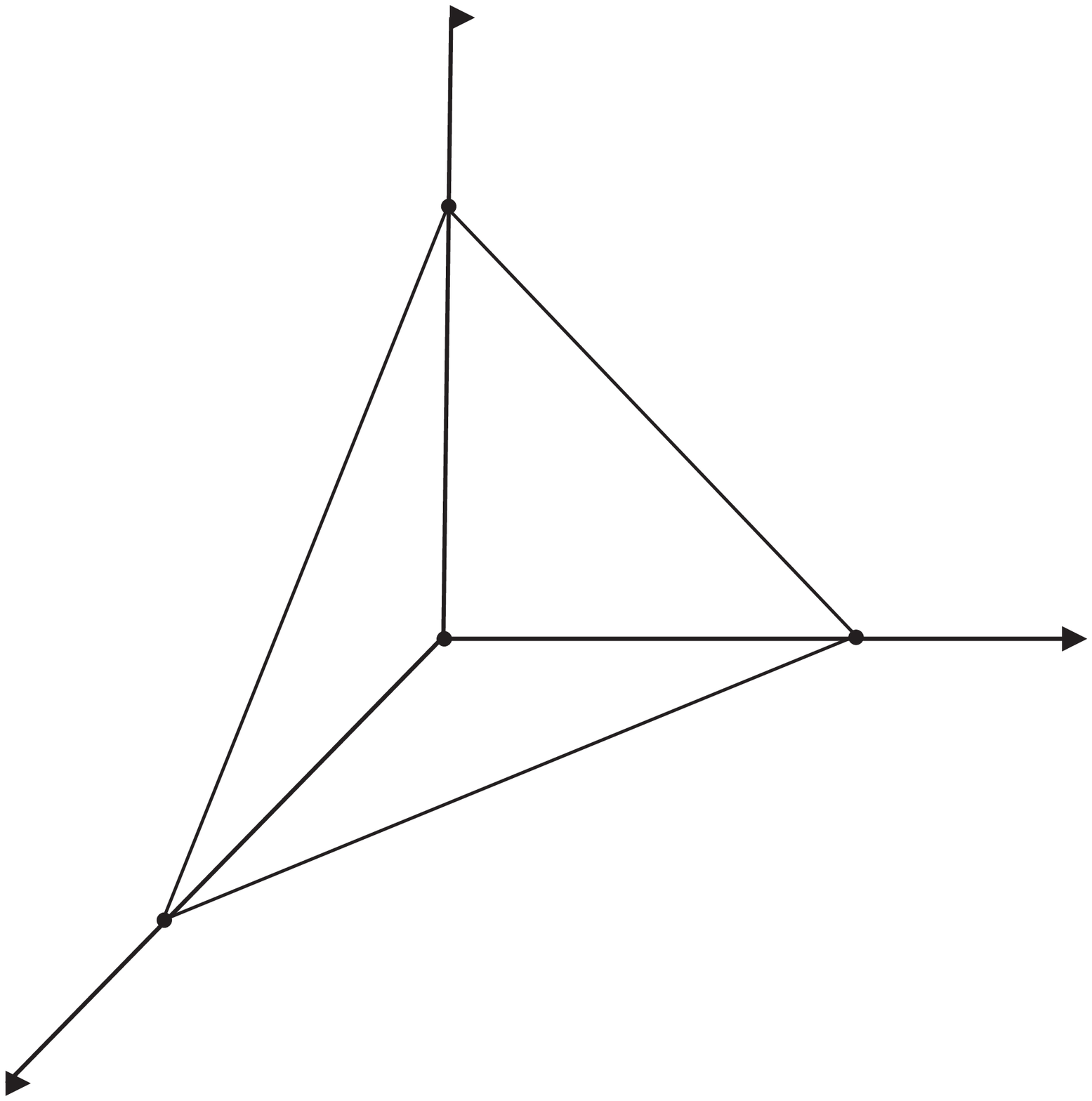}
      \caption{A simplex in $\RR^N$.}
    \end{figure}

Let us say a word about why these facts are true.  Let ${\bf e}_j$ denote
the $j^{\rm th}$ coordinate vector in $\RR^N$ (i.e., the vector with
a 1 in the $j^{\rm th}$ position and 0s in all other slots).  Then a sophisticated
computation with elementary calculus yields that the centroid of the $k$-skeleton
${\cal S}_k$ of the simplex which is the convex hull of $0, {\bf e}_1, {\bf e}_2, \dots, {\bf e}_N$
is
$$
{\cal C}_{k, N} = \frac{1}{N} \cdot \frac{k + (N - k) \sqrt{k + 1}}{(k + 1) + (N - k) \sqrt{k + 1}} ( {\bf e}_1 + {\bf e}_2 +
     \cdots + {\bf e}_N) \, .
$$

From this formula it can immediately be verified that
$$
{\cal S}_0 = {\cal S}_N = \frac{1}{N+1} ( {\bf e}_1 + {\bf e}_2 + \cdots + {\bf e}_N) \, .
$$
It can also be checked that, in dimensions 2 through 12, all the intermediate
skeletons have distinct centroids.  But, in dimension $N = 13$, we observe that
$$
{\cal C}_{3, 13} = {\cal C}_{8, 13} = \frac{23}{13 \cdot 24} ( {\bf e}_1 + {\bf e}_2 +  \cdots + {\bf e}_N) \, .
$$

One may well ask whether dimension $N = 13$ is the only dimension in which there are two intermediate
skeletons with the same centroid.  The answer is ``no''; there are in fact infinitely many such
dimensions (although they are quite sparse---sparser than the prime integers).  One may verify 
this assertion by using the following Diophantine formula.  

\begin{theorem} \sl
Fix a dimension $N \geq 2$.  Consider the simplex {\bf S} as described above.
There are skeletons of dimension $k_1$ and $k_2$, $1 \leq k_1 < k_2 \leq N-1$,
of the simplex {\bf S} which have the same centroid if and only if
$k_1 = a^2 -1$, $k_2 = b^2 - 1$ (for positive integers $a$ and $b$) and, in
addition,
$$
N = (b^2 + ab + a^2) - (b + a) - 1 \, .   \eqno (\star)
$$
\end{theorem}

Obviously this theorem gives us a tool for finding dimensions in which the simplex {\bf S}
has two intermediate skeletons with the same centroid.  The following table gives
some values of the dimension, and of the intermediate dimensions of skeletons which
have the same centroid.  Of course this data may be confirmed by direct calculation
with the formula $(\star)$.  We stress that there are in fact infinitely many
dimensions in which this phenomenon occurs.   The proof of this statement
is a nontrivial exercise in elementary number theory (see [KRMP]).

$$
\begin{array}{|c|c|c|c|}  \hline
\hbox{\bf value of \boldmath $N$} & \hbox{\bf value of
\boldmath $k_1$} & \hbox{\bf value of \boldmath $k_2$}
& \hbox{\bf approx.\ coord.\ of centroid} \\  \hline \hline
13 & 3 & 8 & 0.0737179487 \\
21 & 3 & 15 & 0.0464285714 \\
29 & 8 & 15 & 0.0340038314 \\
31 & 3 & 24 & 0.0317204301 \\
40 & 8 & 24 & 0.0247619047 \\
43 & 3 & 35 & 0.0229789590 \\
51 & 15 & 24 & 0.0194852941 \\
53 & 8 & 35 & 0.0187368973  \\
57 & 3 & 48 & 0.0173872180  \\
65 & 15 & 35 & 0.0153133903 \\  \hline
\end{array}
$$

We conclude this discussion by recording the fact that it is {\it impossible} in 
any dimension for there to be three intermediate skeletons with the same centroid.

\begin{proposition} \sl
For no dimension $N$ can there exists 3 distinct number $1 \leq k_1 < k_2 < k_3 \leq N-1$ such
that the centroids ${\cal C}_{k_1, N}$, ${\cal C}_{k_2, N}$, ${\cal C}_{k_3, N}$ for
the simplex {\bf S} coincide.
\end{proposition}
{\bf Proof:}  We let
$$
Q(a,b) = (b^2 + ab + a^2) - (b + a) - 1 \, .
$$
It suffices for us to show that there do not exist natural numbers $a < b < c$ such that
$Q(a, b) = Q(a, c)$.  Seeking a contradiction, we suppose that such a triple does
indeed exist.

Then
$$
b^2 + ab - b = c^2 + ac - c 
$$
or
$$
b^2 + (a - 1)b = c^2 + (a - 1)c \, .
$$
Since $a \geq 1$, the function $b \mapsto b^2 + (a - 1)b$ is strictly increasing, which yields
a contradiction.
\endpf
\medskip \\

The exploration of centroids for simplices of high dimension is a new venue of exploration.
There are many new phenomena, and more to be discovered.  See [KRMP]  for more results
along these lines.  The reference [ZON] is also of interest.

\newpage

\noindent {\Large \sc References}
\smallskip \\

\begin{enumerate}

\item[{\bf [COT]}]  R. W. Cottle, Quartic barriers, {\it Computational
Optimization and Applications} 12(1999), 81--105.

\item[{\bf [KRA1]}]  S. G. Krantz, A Matter of gravity, {\it Amer.\ Math.\ Monthly} 110(2003), 465--481.

\item[{\bf [KRMP]}]  S. G. Krantz, J. E. McCarthy, and H. R. Parks, Geometric
characterizations of centroids of simplices, {\it Journal of
Mathematical Analysis and Applications} 316(2006), 87--109.

\item[{\bf [ZEE]}]  E. C. Zeeman, {\it Catastrophe Theory.  Selected Papers, 1972--1977},
Addison-Wesley, Reading, MA, 1977.

\item[{\bf [ZON]}]  C. Zong, {\it Strange Phenomena in Convex and Discrete
Geometry}, Springer-Verlag, New York, 1996.

\end{enumerate}
\vspace*{.25in}

\noindent {\bf STEVEN G. KRANTZ} received his B.A. degree from
the University of California at Santa Cruz in 1971. He earned
the Ph.D. from Princeton University in 1974. He has taught at
UCLA, Princeton University, Penn State, and Washington
University in St.\ Louis. Krantz is the holder of the UCLA
Alumni Foundation Distinguished Teaching Award, the Chauvenet
Prize, and the Beckenbach Book Prize. He is the author of 150
papers and 50 books. His research interests include complex
analysis, real analysis, harmonic analysis, and partial
differential equations. Krantz is currently the Deputy
Director of the American Institute of Mathematics.
\medskip \\
{\it American Institute of Mathematics, 360 Portage Avenue, \\
Palo Alto, CA 94306} 
\medskip \\
{\it skrantz@aimath.org}

\end{document}